
\input amstex
\documentstyle{amsppt}
\magnification=\magstep1
\vsize=22.7 truecm
\hsize=16 truecm
 at 14.4truept
\overfullrule=0pt
\def\list#1#2{\noindent\hangindent 3em\hangafter1\hbox to
3em{\hfil#1
              \quad}#2}

\def \g {\frak g}
\def \Ggs {\tilde G\times\tilde\frak g^*}
\def \tG {\tilde G}
\def \tg {\tilde\frak g}
\def \cinf {C^\infty}
\def \ad {{\roman{ad}}}

\def \dt {\left.{d\over dt}\right\vert_{t=0}}
\def \TGk {T^*\tilde G_k}

\def \d {{\roman d}}

\document

\centerline{\bf The Ring of Invariants for Smooth
Completions of  Kac-Moody Lie Algebras}
\smallskip
\bigskip

\centerline{Ian Marshall \footnote{Supported in
part by  CONACyT grant $\#3189-E9307$.}
{}~~~~ and Tudor S. Ratiu
\footnote{Supported in
part by
NSF grant DMS-91-22708 and the Miller Institute for Basic Research at UC
Berkekey.}}
\bigskip

\list{*}{Department of Mathematics, The
University of Leeds, Leeds LS2 9JT, UK}

\list{}{{\it{and}} CIMAT, Apdo.Postal 402, Guanajuato, Gto.36000, Mexico}

\list{**}{Department of Mathematics, University of
California, Santa
Cruz, CA  95064, USA}

\list{} {{\it{and}} IHES, 35, route de Chartres, 91440 Bures-sur-Yvette,
France}
\bigskip

\bigskip
\noindent{\bf Abstract.}~~It is proved that the ring of
invariants of the standard smooth completion of a
Kac-Moody Lie algebra is functionally generated by two elements: the
coefficient
of the center and the Killing form.
\vfil\eject

\noindent{\bf \S 1.~~Introduction}
\smallskip
\noindent
In this paper we address a very specific question:

\medskip
\noindent
What are the smooth ad$^*$-invariant functions for the smooth completion
$\bar\g$ of the full affine Lie algebra
based on a finite dimensional
Lie algebra $\g$?

\medskip
Given is a finite dimensional complex orthogonal Lie algebra $\g$:
that is to say $\g$ admits an invariant symmetric non-degenerate bilinear form.
It gives rise
in the usual fashion
to the smooth completion $\bar \g$ of the affine Lie algebra: $\bar\g$ is
the semidirect product of the complex line $\Bbb C$ with the
standard central extension $\hat\g$ of the loop algebra
$\tg = \cinf(S^1, \g)$.
We will show that the ring of coadjoint-invariant
smooth functions on the dual
$\bar \g^*$ of $\bar\g$ is functionally generated by two elements: the
quadratic form defined by
the naturally induced symmetric bilinear form on $\bar \g$ and the map
associating to an element of $\bar \g^*$ the coefficient of its center.
It is worth
emphasising that when we say generated, we mean generated over $\Bbb C$,
which means that for any $\g$ there are essentially only two complex--valued
coadjoint invariants on $\bar\g$.

A similar result was obtained from an algebraic point of view by Chari and
Ilangovan [1984]. They show that exactly the same two elements polynomially
generate the center of a carefully chosen formal completion of the universal
enveloping algebra of a contragredient Lie algebra defined by a given Cartan
matrix. The present result can be viewed therefore as an extension of
their algebraic theorem to the category of smooth completions of affine
Lie algebras. Other important references are Kac [1984] and Kac and
Peterson [1983].

Our proof is very different to that of Chari and Ilangovan,
being geometric instead of algebraic.
It is based on the idea introduced by Weinstein [1983], of a
dual pair of Poisson maps defined on a symplectic manifold.
 From this geometric point of view, the present
paper can be interpreted as a result on dual pairs for the cotangent
bundle of a loop group.
The natural temptation, based on experience gained from the
finite dimensional case, is to mould the proof in such a way as to be able to
derive the result from standard facts about dual pairs. Unfortunately, as is
mostly the case when working with infinite dimensional objects, such a proof
remains a formal exercise, since the theorems one needs to quote do not in
most cases have infinite dimensional generalizations, and even if they
do, their proofs require the use of technical machinery from the
theory of Fr\'echet manifolds. We have chosen a different way
to deal with these problems: every technical fact for the case at hand
is proved directly, without any appeal to infinite dimensional symplectic or
Poisson geometry. It is important to mention the idea of dual pairs as this
was crucial for
leading us to the result. The interested reader might like to
try and reconstruct our
result from such a point of view, comparing with Marshall [1994],
always bearing in
mind that a demand for total rigour will have to be laid on one side to do
this.

An important motivation for constructing coadjoint invariants on $\bar\g^*$
comes from the field of integrable systems. Invariant functions
on the dual of any
Lie algebra can be used to construct commuting flows.
 From this point of view, the result of the present work is
disappointing for it rules out the search for interesting
integrable systems based on
$\bar\g$ (at least by the standard approach).
However, on
both $\g^*$ and $\hat\g^*$ (these are the duals respectively to $\g$ and to the
central extension of the loop algebra based on $\g$),
there are a lot of invariant functions which have been put to good use in the
works on integrable  systems.

\bigskip

\noindent{\bf \S 2.~~Statement of the theorem}

\smallskip

Let $\g$ be a finite dimensional complex Lie algebra endowed with a
nondegenerate
ad--invariant inner product $\bigl(\ ,\ \bigr):\g\times\g\rightarrow\g$,
that
is,
$$
\bigl(X,[Y,Z]\bigr)=\bigl([X,Y],Z\bigr)\quad\forall X,Y,Z\in\g.
\eqno(2.1)
$$
The full affine Lie algebra $\bar\g$ based on $\g$ is obtained by the
following
construction (see Kac [1985], Pressley and Segal [1986]). Consider the
central extension $\hat\g = \tg\oplus\Bbb C$ of the loop algebra
$\tg =\cinf(S^1,\g)$
by
$\Bbb C$, having the Lie bracket $$
[(X,a),(Y,b)](x)=\left([X(x),Y(x)],\,
\int_0^{2\pi}\bigl(X(y),Y'(y)\bigr)dy\right),
\eqno(2.2)
$$
where  the prime ``$\ '\ $'' denotes derivation. Let the abelian Lie
algebra $\Bbb C$ act on
$\hat\g$ by
$$
z\cdot(X,a)=(zX',0).
\eqno(2.3)
$$
and form the semidirect product $\bar\g=\Bbb C\oplus\hat\g$ with Lie bracket
$$
[(z_1,X,a)\,,(z_2,Y,b)]=
\left(0,\,[X,Y] + z_1Y' - z_2X',\, \int_0^{2\pi}\bigl(X(y),
Y'(y)\bigr)dy\right).\eqno(2.4)
$$

If we identify $\bar\g^*$ with $\bar\g$ by
$$
<(\alpha,\xi,e),\,(z,X,a)> =
Re\left(\alpha z + \int_0^{2\pi}\bigl(\xi(x),X(x)\bigr)dx + ea\right)
\eqno(2.5)
$$
then we obtain
$$
\ad^*_{(z,X,a)}(\alpha,\xi,e)
= \left(\int_0^{2\pi}\bigl(\xi(x),X'(x)\bigr)dx,\, [X,\xi] + z\xi' +
eX',\,0\right), \eqno(2.6)
$$
where $\ad^*$ denotes minus the dual map of $\ad$. We will denote by the
same
symbol $<,>$ the pairing identifying $\tg^*$ with $\tg$, i.e., take
in (2.5) $\alpha = e = z = a = 0$.

If $F \in \cinf(\bar\g^*)$,
we
will denote by  $\nabla_{(\alpha, \xi, e)}F \in \bar\g$
its functional derivative,
defined to be the vector representative of the  Fr\'echet
derivative of $F$ at $(\alpha, \xi, e)$ relative to the
pairing (2.5).That is
$$
\eqalignno{
f\bigl((\alpha,\,\xi,\,e) + t(\beta,\,\eta,&\, e')\bigr)\ =\
f(\alpha,\,\xi,\,e) + t<(\beta,\,\eta,\,e'),\nabla_{(\alpha,\,\xi,\,e)}
f>\cr
&+ O(t^2)\ \ \roman{as}\ \ t\rightarrow 0,\ t\in\Bbb R\quad
\forall\ (\alpha,\,\xi,\,e), (\beta,\,\eta,\,e')\in\bar\g^*.
&(2.7)\cr}
$$
Our goal is to find the set $I(\bar\g^*)$
($I$ stands for {\it invariant})
of
all functions $F \in \cinf(\bar\g^*)$ which are $\ad^*$--invariant, that is,
satisfying
$$ <(\alpha,\xi,e),\, [(z,X,a), \nabla_{(\alpha, \xi, e)}F ]> = 0 \quad
\forall (z,X,a)\in\bar\g,\forall(\alpha,\xi,e)\in\bar\g^*.
\eqno (2.8) $$
\smallskip
\noindent We will prove the following theorem:
\smallskip
\noindent {\bf Theorem} {\it Any} $\ad^*${\it--invariant
function on $\bar\g^*$ is
of the form $\Cal F(\kappa, \pi)$, where $\kappa:
\bar\g^*\rightarrow\Bbb C$
is
given by $$
\kappa(\alpha,\xi,e)=e\alpha -
{1\over2}\int_0^{2\pi}\bigl(\xi(x),\xi(x)\bigr)dx,
\eqno (2.9)$$
$\pi: \bar\g^* \rightarrow \Bbb C$ is the projection onto the center,
and $\Cal F\in\cinf(\Bbb C^2)$.}

\smallskip
Harnad and Kupershmidt [1991] is also relevant as it contains
many of the formulae for the momentum maps in section 3.
\bigskip

\noindent{\bf \S 3.~~Proof of the theorem}
\smallskip

Suppose that $G$ is a connected and simply connected Lie group whose
Lie algebra is $\g$, regarded as a real Lie algebra.
Let $\tG =
\cinf(S^1,G)$ be the loop group of $G$ and let $\tg = \cinf(S^1,\g)$ be
the corresponding Lie algebra. The loop group $\tG$ is connected. On the
cotangent bundle
$T^*\tG$,  which we will henceforth identify with $\Ggs$ by left
translations we
define the
family of non--canonical symplectic structures by
$$
\omega_k =  \omega_{CAN} + \pi^*\Omega_k\,,\quad\quad k\in\Bbb C.
\eqno(3.1)
$$
Here $\omega_{CAN}$ is the canonical 2-form and  $\pi: T^*\tilde
G\rightarrow\tG$ is the canonical projection to the base of the bundle.
$\pi^*$ denotes the pull-back defined by $\pi$. $\Omega_k$ is
the closed 2-form on $\tG$ given by
$$
\Omega_k(g)(V,W) = Re\int_0^{2\pi}\bigl(X(x),kY'(x)\bigr)(x)dx = <X, kY'>
\eqno(3.2)
$$
where $V,W\in T_g\tG$ are vectors tangent to the curves
$$
t\mapsto\ \ ge^{tX},\ \text{and} \quad t\mapsto\   \ ge^{tY}
\eqno(3.3)
$$
respectively.
We will refer to the symplectic manifold
$(T^*\tG,\ \omega_k)$ as $\TGk$.

The left and right actions of $\tG$ on itself lift to the following
symplectic
actions on $\TGk$:
$$
\tG\times \TGk\rightarrow \TGk,\qquad
(h,(g,\mu))\mapsto L_h(g,\mu) = (hg,\ \mu)
\eqno(3.4)
$$
and
$$
\tG\times \TGk\rightarrow \TGk,\qquad
(h,(g,\mu))\mapsto R_{h^{-1}}(g,\mu) = (gh^{-1},\ h\mu h^{-1} + kh'h^{-1}).
\eqno(3.5)
$$
The corresponding momentum maps are given respectively by
$$
J^L(g,\mu) = g\mu g^{-1} + kg'g^{-1}
\eqno(3.6)
$$
and
$$
J^R(g,\mu) = -\mu.
\eqno(3.7)
$$
See Section 4. In the above formulae $h\mu h^{-1}$
means Ad$^*_h\mu$ and similarly for
$ g\mu g^{-1}$, where Ad$^*$ is the coadjoint action of
$G$ on $\g^*$ (that is
Ad$^*$  is the dual of the inverse of Ad). Below, we
prefer to reserve the notations
Ad$^*$ and ad$^*$ for actions on $\bar\g^*$.
In order not to clutter the
flow
of the proof, we will prove (3.6) and (3.7)
separately in Section 4.

Let us now observe that $J^L$ and $J^R$ are independent on $\TGk$. To see
this it is sufficient to check that for $\xi,\eta\in\tg$ the
function $\Phi=<\xi,J^L>+
<J^R,\eta>$ cannot have zero exterior derivative
in any open set unless $\xi=0=\eta$.
If $V\in T_{(g,\mu)}\TGk$ is tangent to the curve
$t\mapsto(ge^{tX},\mu+t\nu)$, we get
$$
\eqalignno{
(V\lrcorner \d\Phi)(g \mu) &= \dt\Phi(ge^{tX},\mu+t\nu)\cr
    &= <g^{-1}\xi g, \nu+[X,\mu]+kX'> - <\eta,\nu>,
&(3.8)\cr}
$$
\noindent
from where the assertion follows.

Let us also compute the Poisson brackets. For any $X,Y\in\tg$ we have
$$
\{<J^L,X>,<J^L,Y>\} = <J^L,[X,Y]> + <X,kY'>,
\eqno(3.9a)
$$
$$
\{<J^R,X>,<J^R,Y>\} = <J^R,[X,Y]> - <X,kY'>,
\eqno(3.9b)
$$
$$
\{<J^L,X>,<J^R,Y>\} = 0.
\eqno(3.9c)
$$
It follows that all Poisson brackets of the form
$\{{J^L}^*\phi,{J^R}^*\psi\}$
vanish, where ${J^L}^*$ and ${J^R}^*$ denote the respective pull-back
mappings.
Note also that it follows from (3.9a) and from (3.9b) that both $J^L$ and
$J^R$
are not infinitesimally equivariant. See the last paragraph of section 4.

We will denote by $\ell(X)$ and $r(X)$ the
infinitesimal generators of the left
and right actions (3.4) and (3.5); they are vector
fields on $\TGk$. Thus, by
definition of the momentum map,
$$ \ell(X) = \Bbb X_{<J^L,X>}\quad{\roman{and}}\quad
r(X) = \Bbb X_{<J^R,X>},
\eqno(3.10)
$$
where $\Bbb X_H$ denotes the Hamiltonian vector field
corresponding to the function $H$.

Suppose now that $\Phi\in\cinf(\TGk)$ and that $\Phi$ can be factored
through
$J^L\times J^R$. If $\Phi$ is invariant
under the left action $L$ given by (3.4),
then $\Phi$ is a function of $J^R$ only. That is, if $\Phi = f(J^L,J^R)$
then
$$
\Phi\circ L_h = \Phi\quad\forall h\in\tG\quad
\Leftrightarrow\quad
{{\partial f}\over{\partial a}}(a,b) = 0 \quad\forall a,b\in\tg^*.
\eqno(3.11)
$$

We define the infinitesimal action $v$ of
$\Bbb C$, considered to be a real abelian
Lie algebra, on $\TGk$:
$$
v:\Bbb C\rightarrow vect(\TGk),
\eqno(3.12)
$$
by declaring $v(z)(g,\mu)\in T_{(g,\mu)}\TGk$
to be the vector tangent at $t=0$ to
the curve
$$
t \mapsto (ge^{tzg^{-1}g'},\mu+tz\mu').
\eqno(3.13)
$$
Clearly $v$ defines an infinitesimal
action, i.e.  $[v(z_1),v(z_2)]=0$ for any
$z_1,z_2\in\Bbb C$. It is also
an infinitesimally symplectic action, since for any
$F,\ H\in\cinf(\TGk)$ we have
$$
\{v(z)\lrcorner \d F,\,H\} + \{F,\,v(z)
\lrcorner \d H\}
=
v(z)\lrcorner\d\{F,H\},
\eqno(3.14)
$$
where $\lrcorner$ denotes the interior derivative
operation (contraction on the
first index) of a vector field on a form.  The momentum map
for the action given by (3.13) is
computed in Section 4 and has the expression
$$
\Cal J(g,\mu) = \int_0^{2\pi}(g^{-1}g',\,\mu) + {1\over2}
\int_0^{2\pi}
(kg^{-1}g',g^{-1}g').
\eqno(3.15)
$$
The map $\Cal J$ is clearly infinitesimally equivariant.
If $k\neq0$, (3.15) can be
rewritten in the following way:
$$
\Cal J = {1\over2k}\int_0^{2\pi}\left((J^L,J^L)-(J^R,J^R)\right).
\eqno(3.16)
$$
Note that the algebraic dual to $\Bbb C$ has been identified with
$\Bbb C$ itself via the real pairing that associates to $(z_1, z_2)$
the real part of the product $z_1z_2$.

The $\Bbb C$
action $v$ and the left action $L$ of $\tG$ on $\TGk$ are compatible
in the sense that
$$ [v(z), \ell(X)] = - \ell(z\cdot X)
\eqno(3.17)
$$
where $z\cdot X = zX'$. Therefore the semi--direct product
$\Bbb C\ltimes\tg$ acts infinitesimally symplectically
on $\TGk$ and this action has the momentum map
$$
\Cal J\oplus J^L\ :\ \TGk\ \rightarrow\ \Bbb C\times\tg^*.
\eqno(3.18)
$$
This is a {\it{non-equivariant}} momentum map from $\TGk$
to the dual of the semi--direct product
$\Bbb C\ltimes\tg$.
As recalled in Section 4, the cocycle
term needed to define the canonical
equivariant extension of this momentum
map is found by the following
computation:
$$
\eqalignno{
\{<&(\Cal J, J^L),(z,X)>,\
<(\Cal J, J^L),(\zeta,Y)>\}\cr
&=
<J^L,\ [X,Y] + zY' - \zeta X'> + <X,kY'>\cr
&=<(\Cal J, J^L),
[(z, X), (\zeta, Y)]> + <X,kY'>\cr
&=<(\Cal J,J^L,k),
[(z, X,a), (\zeta, Y,b)]>.
&(3.19)\cr
}
$$
where the Lie bracket in the final expression is that given by (2.4) and $a$
and
 $b$ are any complex numbers.
In other words, the map

$$
\frak S = (\Cal J, J^L,k):\TGk\ \rightarrow\ \bar\g^*
\eqno(3.20)
$$
is infinitesimally equivariant, i.e for any
$\varphi\in \cinf(\bar\g^*)$,
and for any $(z, X, \sigma)\in\bar\g$,
$$
(v(z)+\ell(X))\lrcorner{\roman d}
(\varphi\circ\frak S)
=
<\frak S,\,[\nabla_\frak S \varphi,\, (z,X, \sigma)]>.
\eqno(3.21)
$$
\smallskip

Suppose now that $F\in I(\bar\g^*)$. By (2.8) we have
$$
<(\alpha,\xi,e),\ [(z,X,\sigma),\,
\nabla_{(\alpha,\xi,e)}F]> = 0
\quad\forall(z,X,\sigma)\in\bar\g,\
\forall(\alpha,\xi,e)\in\bar\frak
g^*.
\eqno(3.22)
$$
Letting $\frak S(g,\mu) = (\alpha, \xi,e)$ and
combining relations (3.21) and (3.22) for $z=0$, we get
$$
\ell(X)\lrcorner{\roman d}(F\circ\frak S)=0\quad\forall X\in\tg.
\eqno(3.23)
$$

\smallskip
\noindent
Thus we have shown that $F\circ\frak S\in \cinf(\TGk)$
is invariant with respect to
the infinitesimal left action of $\tg$ on $\TGk$. Moreover, from (3.16),
$F\circ\frak S$ factors through $J^L\times J^R$. By (3.11) it
follows then that $F\circ\frak S$
must be of the form
$$
F\circ\frak S = \psi\circ J^R
\eqno(3.24)
$$
for some $\psi\in \cinf(\tg^*)$, i.e. we have
$$
\eqalignno{
\psi\circ J^R  &= F(\Cal J,J^L,k)\cr
&=F\bigl( {1\over2k}\int_0^{2\pi}\bigl((J^L,J^L)\bigr)-
(J^R,J^R)\bigr),
J^L,k \bigr).
&(3.25)\cr
}
$$
i.e. $F$ is a function of $\int_0^{2\pi}(J^R,J^R)$.
Denoting $J^L(g, \mu) = \xi,\, {\Cal J}(g, \mu) = \alpha$
and observing that
$$
-{1\over2}\int_0^{2\pi}(J^R,J^R)(g, \mu) = k\alpha - {1\over2}\int_0^{2\pi}
(\xi, \xi),
\eqno(3.26)
$$
it follows that for each fixed $k$, $F$ has the form
$$
F = \Cal G\circ\kappa
\eqno(3.27)
$$
for some $\Cal G\in \cinf(\Bbb C)$, where
$$\kappa(\alpha,\xi,k)=k\alpha -
{1\over2}\int_0^{2\pi}\bigl(\xi(x),\xi(x)\bigr)dx.
\eqno (3.28)
$$
Functions depending only on $k$ are clearly invariant since $k$ is the
coordinate
of the center. Therefore, the above result holds for each $k$ and we
conclude
that
$$F = \Cal F (\kappa, \pi),
\eqno(3.29)$$
where $\pi: \bar\g^* \rightarrow \Bbb C$ is the projection
onto the center, $\kappa$ is given by (3.28)
and $\Cal F\in\cinf(\Bbb C^2)$.

\bigskip

\noindent{\bf \S 4.~~The momentum maps $J^L,\ J^R,\ \Cal J$}
\smallskip

For an introduction to symplectic manifolds
the reader is
referred to one of the many standard books on
the subject, for example
Abraham and Marsden [1978] or Arnold [1978].
Let us begin with a closer look at the symplectic manifold
$\TGk$.

Let $\varphi, \psi \in \cinf (\Ggs)$ and denote
by $\varphi^\mu, \psi^\mu \in \cinf (\tG), \varphi^g, \psi^g \in
\cinf (\tg )$ the partial functions $\varphi^\mu(g) = \varphi^g(\mu)
=\varphi (g, \mu), \psi^\mu(g) = \psi^g(\mu) = \psi (g, \mu)$. The
Poisson bracket is given by
$$
\eqalignno{
\{\varphi, \psi\}(g,\mu) = &\dt\bigl(\varphi^\mu(ge^{t\nabla\psi^g(\mu)})
- \psi^\mu(ge^{t\nabla\varphi^g(\mu)})\bigr)\cr
&-<\mu, [\nabla \varphi^g(\mu), \nabla \psi^g(\mu)]> -
<k\nabla \varphi^g(\mu), (\nabla \psi^g(\mu))'>
&(4.1)\cr
}
$$
where $\nabla \varphi^g(\mu), \nabla \psi^g(\mu) \in \tg$
are the gradients of $\varphi^g$ and $\psi^g$ relative to the
pairing $<,>$. (The reader might find it helpful to note that the rule for
computing the gradient of a function relative to a pairing was given explicitly
in (2.7) for the pairing given by (2.5).)

It is straightforward to verify that $L_h$ and
$R_{h^{-1}}$ define symplectic
actions by  checking,
$$
\{\varphi\circ L_h,\psi\circ L_h\} = \{\varphi,\psi\}
\circ L_h\quad
\forall\varphi,\ \psi\in \cinf(\Ggs),\ {\roman{for}}\
{\roman{any}}\ h\in\tilde G
\eqno(4.2)
$$
and similarly for $R_{h^{-1}}$.

The values at $(g, \mu)$ of the infinitesimal
generators of these two actions,
$\ell(X)$ and $r(X)$, are the
derivatives with respect to $t$ at $t=0$ to
the curves in $\Ggs$ given by
$$t \mapsto
(e^{tX}g,\mu) \quad \text{and} \quad t
\mapsto (ge^{-tX},e^{tX}\mu
e^{-tX}+tkX').
\eqno(4.3)
$$
By the definition of the momentum map, $J^L$ and
$J^R$ are given as the solutions to
the equations
$$
\Bbb X_{<J^L,X>} = \ell(X)\quad\quad
\Bbb X_{<J^R,X>} = r(X),
\eqno(4.4)
$$
where $\Bbb X_H\in vect(\Ggs)$ is the Hamiltonian
vector field corresponding
to the function $H$, i.e.
$$
\Bbb X_H\lrcorner\d K = \{K,H\},
\quad\forall K\in \cinf(\Ggs).
\eqno(4.5)
$$
Let us verify (4.4) for the map $J^L$.
If $\varphi\in\cinf(\tG\times\tg^*)$, we have

$$
\eqalignno{
(\ell(X)\lrcorner\d \varphi)(g,\mu)\ &=\
\dt \varphi(e^{tX}g,\mu)\cr
&= \{\varphi,<J^L, X>\}(g,\mu).
&(4.6)\cr}
$$
The last equality in (4.6) is obtained
by using (4.1) and the formulae
$$
\eqalignno{
\nabla <J^L, X>^g(\mu) &= g^{-1}Xg,\cr
\dt<J^L, X>^{\mu}(ge^{tY})&=
<[Y, \mu] + kY', g^{-1}Xg>\ \ {\roman{for}}\ {\roman {any}}\ Y\in\tg.
&(4.7)\cr}
$$
One similarly verifies formula (4.4) for the momentum map
$J^R$ by using (4.1) and
$$
\eqalignno{
\nabla <J^R, X>^g (\mu) &= -X,\cr
\dt<J^R, X>^{\mu}(ge^{tY}) &= 0\ \ {\roman{for}}\ {\roman{any}}\ Y\in\tg.
&(4.8)\cr}
$$

\smallskip
Next, we turn to the computation of the
momentum map for the $\Bbb C$ action $v$.
For $\varphi\in\cinf(\tG\times\tg^*)$, we have
$$
\eqalignno{
(v(z)\lrcorner\d\varphi)(g,\mu)\ &=\
\dt\varphi(ge^{tzg^{-1}g'},\mu+tz\mu')\cr
&= \{\varphi,<\Cal J,z>\}(g,\mu),
&(4.9)\cr}
$$
which again follows from (4.1) and the formulae

$$
\nabla <\Cal J, z>^g(\mu) = zg^{-1}g',
\eqno(4.10a)
$$
$$
\eqalignno{
\dt<\Cal J, z>^\mu(ge^{tY}) &=\
<[g^{-1}g', Y] + Y', z\mu> +\cr
<[g^{-1}&g', Y] + Y', zkg^{-1}g'>\ \ {\roman{for}}\ {\roman{any}}\ Y\in\tg.
&(4.10b)\cr}
$$

\medskip

In general, if $P:M\rightarrow\frak a^*$ is a
momentum map for the infinitesimal symplectic
action of a Lie algebra $\frak a$ on a connected
symplectic manifold $M$, we
have

$$
[\Bbb X_{<P,A>},\Bbb X_{<P,B>}] = \Bbb X_{<P,[AB]>}
\quad\forall A,B\in\frak a,
\eqno(4.11)
$$
but the left hand side of $(4.11)$ is
$\Bbb X_{\{<P,A>,<P,B>\}}$ and it follows
that
$$
\{<P,A>,<P,B>\} = <P,[A,B]> + \, \Sigma(A,B),
\eqno(4.12)
$$
where $\Sigma(A,B)$ is a constant function on $M$.
It can be checked that the constant term in
(4.12) corresponds to a coadjoint
one-cocycle on $\frak a$ and hence although
$P$ may not be equivariant
(equivariance would be the case when the
constant is zero), an equivariant  map
$\hat P$ can be constructed from $P$, having
its image in the central extension
$\hat\frak a = \frak a\oplus\Bbb R$ of
$\frak a$. This is the reason for the
computation in (3.19) and for the remark following
the formulae in (3.9).

\bigskip

\centerline{\bf REFERENCES}
\medskip

\noindent {R. Abraham and J.Marsden [1978] {\sl Foundations of Mechanics},
Benjamin/Cummings Publishing Co. Reading, Mass.}

\noindent {V.I. Arnold [1978] {\sl Mathematical
Methods of Classical Mechanics}, Springer--Verlag.
New York, Heidelberg}

\noindent {V. Chari and S. Ilangovan [1984] On the Harishchandra
homomorphism for infinite dimensional Lie algebras,
{\sl Journal of Algebra \/}, {\bf 90}, 476--490.}

\noindent
V.G. Kac and Peterson [1983]
Infinite flag varieties and conjugacy theorems,
{\sl Proceedings of the National Academy of Sciences\/},
{\bf 80}, 1778--1782.

\noindent
V.G. Kac [1984]
Laplace operators of infinite-dimensional Lie algebras and
theta functions,
{\sl Proceedings of the National Academy of Sciences\/},
{\bf 81}, 645--647.

\noindent {V.G. Kac [1985] {\sl Infinite Dimensional Lie Algebras\/},
Cambridge University Press, Cambridge, New York}

\noindent {I. Marshall [1994]
 Modified Systems found by symmetry reduction on the cotangent bundle of a
loop
group, {\sl J. Geom. and Phys.},to appear.}

\noindent {J.Harnad \& B.A,Kupershmidt [1991] Symplectic geometries on
$T^*\tilde G$, Hamiltonian group actions and integrable systems. {\sl
Preprint}}.

\noindent {A.Pressley and G. Segal [1986] {\sl
Loop Groups\/}, Oxford University Press, Oxford, New York}

\noindent {A. Weinstein [1983] The local
structure of Poisson manifolds,
{\sl J. Diff. Geom.\/},{\bf 18}, 523--557.}

\enddocument